\documentclass{amsart}
\usepackage{amssymb,amsmath,theorem,euscript}

\input{epsf}

\def\sm{\smallskip}

\begin{document}

\def\wh{\widehat}

 \def\ov{\overline}
\def\wt{\widetilde}
 \newcommand{\rk}{\mathop {\mathrm {rk}}\nolimits}
\newcommand{\Aut}{\mathop {\mathrm {Aut}}\nolimits}
\newcommand{\Out}{\mathop {\mathrm {Out}}\nolimits}
 \newcommand{\tr}{\mathop {\mathrm {tr}}\nolimits}
  \newcommand{\diag}{\mathop {\mathrm {diag}}\nolimits}
  \newcommand{\supp}{\mathop {\mathrm {supp}}\nolimits}
  \newcommand{\indef}{\mathop {\mathrm {indef}}\nolimits}
  \newcommand{\dom}{\mathop {\mathrm {dom}}\nolimits}
  \newcommand{\im}{\mathop {\mathrm {im}}\nolimits}
 
\renewcommand{\Re}{\mathop {\mathrm {Re}}\nolimits}

\def\Br{\mathrm {Br}}

\def\SL{\mathrm {SL}}
\def\Diag{\mathrm {Diag}}
\def\SU{\mathrm {SU}}
\def\GL{\mathrm {GL}}
\def\U{\mathrm U}
\def\OO{\mathrm O}
 \def\Sp{\mathrm {Sp}}
 \def\SO{\mathrm {SO}}
\def\SOS{\mathrm {SO}^*}
 \def\Diff{\mathrm{Diff}}
 \def\Vect{\mathfrak{Vect}}
\def\PGL{\mathrm {PGL}}
\def\PU{\mathrm {PU}}
\def\PSL{\mathrm {PSL}}
\def\Symp{\mathrm{Symp}}
\def\End{\mathrm{End}}
\def\Mor{\mathrm{Mor}}
\def\Aut{\mathrm{Aut}}
 \def\PB{\mathrm{PB}}
 \def\cA{\mathcal A}
\def\cB{\mathcal B}
\def\cC{\mathcal C}
\def\cD{\mathcal D}
\def\cE{\mathcal E}
\def\cF{\mathcal F}
\def\cG{\mathcal G}
\def\cH{\mathcal H}
\def\cJ{\mathcal J}
\def\cI{\mathcal I}
\def\cK{\mathcal K}
 \def\cL{\mathcal L}
\def\cM{\mathcal M}
\def\cN{\mathcal N}
 \def\cO{\mathcal O}
\def\cP{\mathcal P}
\def\cQ{\mathcal Q}
\def\cR{\mathcal R}
\def\cS{\mathcal S}
\def\cT{\mathcal T}
\def\cU{\mathcal U}
\def\cV{\mathcal V}
 \def\cW{\mathcal W}
\def\cX{\mathcal X}
 \def\cY{\mathcal Y}
 \def\cZ{\mathcal Z}
\def\0{{\ov 0}}
 \def\1{{\ov 1}}
 \def\frA{\mathfrak A}
 \def\frB{\mathfrak B}
\def\frC{\mathfrak C}
\def\frD{\mathfrak D}
\def\frE{\mathfrak E}
\def\frF{\mathfrak F}
\def\frG{\mathfrak G}
\def\frH{\mathfrak H}
\def\frI{\mathfrak I}
 \def\frJ{\mathfrak J}
 \def\frK{\mathfrak K}
 \def\frL{\mathfrak L}
\def\frM{\mathfrak M}
 \def\frN{\mathfrak N} \def\frO{\mathfrak O} \def\frP{\mathfrak P} \def\frQ{\mathfrak Q} \def\frR{\mathfrak R}
 \def\frS{\mathfrak S} \def\frT{\mathfrak T} \def\frU{\mathfrak U} \def\frV{\mathfrak V} \def\frW{\mathfrak W}
 \def\frX{\mathfrak X} \def\frY{\mathfrak Y} \def\frZ{\mathfrak Z} \def\fra{\mathfrak a} \def\frb{\mathfrak b}
 \def\frc{\mathfrak c} \def\frd{\mathfrak d} \def\fre{\mathfrak e} \def\frf{\mathfrak f} \def\frg{\mathfrak g}
 \def\frh{\mathfrak h} \def\fri{\mathfrak i} \def\frj{\mathfrak j} \def\frk{\mathfrak k} \def\frl{\mathfrak l}
 \def\frm{\mathfrak m} \def\frn{\mathfrak n} \def\fro{\mathfrak o} \def\frp{\mathfrak p} \def\frq{\mathfrak q}
 \def\frr{\mathfrak r} \def\frs{\mathfrak s} \def\frt{\mathfrak t} \def\fru{\mathfrak u} \def\frv{\mathfrak v}
 \def\frw{\mathfrak w} \def\frx{\mathfrak x} \def\fry{\mathfrak y} \def\frz{\mathfrak z} \def\frsp{\mathfrak{sp}}
 \def\bfa{\mathbf a} \def\bfb{\mathbf b} \def\bfc{\mathbf c} \def\bfd{\mathbf d} \def\bfe{\mathbf e} \def\bff{\mathbf f}
 \def\bfg{\mathbf g} \def\bfh{\mathbf h} \def\bfi{\mathbf i} \def\bfj{\mathbf j} \def\bfk{\mathbf k} \def\bfl{\mathbf l}
 \def\bfm{\mathbf m} \def\bfn{\mathbf n} \def\bfo{\mathbf o} \def\bfp{\mathbf p} \def\bfq{\mathbf q} \def\bfr{\mathbf r}
 \def\bfs{\mathbf s} \def\bft{\mathbf t} \def\bfu{\mathbf u} \def\bfv{\mathbf v} \def\bfw{\mathbf w} \def\bfx{\mathbf x}
 \def\bfy{\mathbf y} \def\bfz{\mathbf z} \def\bfA{\mathbf A} \def\bfB{\mathbf B} \def\bfC{\mathbf C} \def\bfD{\mathbf D}
 \def\bfE{\mathbf E} \def\bfF{\mathbf F} \def\bfG{\mathbf G} \def\bfH{\mathbf H} \def\bfI{\mathbf I} \def\bfJ{\mathbf J}
 \def\bfK{\mathbf K} \def\bfL{\mathbf L} \def\bfM{\mathbf M} \def\bfN{\mathbf N} \def\bfO{\mathbf O} \def\bfP{\mathbf P}
 \def\bfQ{\mathbf Q} \def\bfR{\mathbf R} \def\bfS{\mathbf S} \def\bfT{\mathbf T} \def\bfU{\mathbf U} \def\bfV{\mathbf V}
 \def\bfW{\mathbf W} \def\bfX{\mathbf X} \def\bfY{\mathbf Y} \def\bfZ{\mathbf Z} \def\bfw{\mathbf w}
 \def\R {{\mathbb R }} \def\C {{\mathbb C }} \def\Z{{\mathbb Z}} \def\H{{\mathbb H}} \def\K{{\mathbb K}}
 \def\N{{\mathbb N}} \def\Q{{\mathbb Q}} \def\A{{\mathbb A}} \def\T{\mathbb T} \def\P{\mathbb P} \def\G{\mathbb G}
 \def\bbA{\mathbb A} \def\bbB{\mathbb B} \def\bbD{\mathbb D} \def\bbE{\mathbb E} \def\bbF{\mathbb F} \def\bbG{\mathbb G}
 \def\bbI{\mathbb I} \def\bbJ{\mathbb J} \def\bbK{\mathbb K} \def\bbL{\mathbb L} \def\bbM{\mathbb M} \def\bbN{\mathbb N} \def\bbO{\mathbb O}
 \def\bbP{\mathbb P} \def\bbQ{\mathbb Q} \def\bbS{\mathbb S} \def\bbT{\mathbb T} \def\bbU{\mathbb U} \def\bbV{\mathbb V}
 \def\bbW{\mathbb W} \def\bbX{\mathbb X} \def\bbY{\mathbb Y} \def\kappa{\varkappa} \def\epsilon{\varepsilon}
 \def\phi{\varphi} \def\le{\leqslant} \def\ge{\geqslant}

\def\UU{\bbU}
\def\Mat{\mathrm{Mat}}
\def\tto{\rightrightarrows}

\def\Gr{\mathrm{Gr}}

\def\graph{\mathrm{graph}}

\def\O{\mathrm{O}}

\def\la{\langle}
\def\ra{\rangle}

\def\B{\mathrm B}
\def\Int{\mathrm{Int}}
\def\LGr{\mathrm{LGr}}


\def\I{\mathbb I}
\def\M{\mathbb M}
\def\T{\mathbb T}
\def\S{\mathrm S}

\def\Lat{\mathrm{Lat}}
\def\LLat{\mathrm{LLat}} 
\def\Mod{\mathrm{Mod}}
\def\LMod{\mathrm{LMod}}
\def\Naz{\mathrm{Naz}}
\def\naz{\mathrm{naz}}
\def\bNaz{\mathbf{Naz}}
\def\AMod{\mathrm{AMod}}
\def\ALat{\mathrm{ALat}}
\def\MAT{\mathrm{MAT}}
\def\Mar{\mathrm{Mar}}

\def\Ver{\mathrm{Vert}}
\def\Bd{\mathrm{Bd}}
\def\We{\mathrm{We}}
\def\Heis{\mathrm{Heis}}
\def\Pol{\mathrm{Pol}}
\def\Ams{\mathrm{Ams}}
\def\Herm{\mathrm{Herm}}

 \def\kos{/\!\!/}
  \newcommand{\sgn}{\mathop {\mathrm {sgn}}\nolimits}
  
  \renewcommand{\Re}{\mathop {\mathrm {Re}}\nolimits}
\renewcommand{\Im}{\mathop {\mathrm {Im}}\nolimits}

\def\Gr{\mathrm{Gr}}

\begin{center}
 \Large\bf

After Plancherel formula

\smallskip

\large \sc

\sm

Yury Neretin%
 \footnote{Supported by the grant FWF, P28421.}

\end{center}

\smallskip

{\small We discuss two topics related to Fourier transforms
on Lie groups and on homogeneous spaces:
the operational calculus and the Gelfand--Gindikin problem (program) about separation of non-uniform spectra.
Our purpose is to indicate some  non-solved problems of non-commutative harmonic analysis
that definitely are solvable.
This is a sketch of my talks on VI School "Geometry and Physics", Bia\l{}owie\.{z}a, Poland, June 2017.}

\smallskip

{\bf 1. Abstract Plancherel theorem for groups.} See, e.g., \cite{Fol}.
Let $G$ be a type
I locally compact group with a two-side invariant Haar  measure $dg$. Denote by $\wh G$ the set of all irreducible unitary representations of
$G$ (defined up to a unitary equivalence%
\footnote{For a formal definition of type I groups see. e.g., \cite{Fol}, Sect. 7.2.
 Connected semisimple Lie groups, connected nilpotent Lie groups,
classical $p$-adic groups have type I.
This condition implies  a presence of the standard Borel structure on $\wh G$
and a uniqueness of a decomposition
of any unitary representation of $G$ into a direct integral of irreducible representations.}).
For $\rho\in \wh G$ denote by $H_\rho$ the space of the representation $\rho$. For $\rho\in \wh G$ and $f\in L^1(G)$ 
we define the following operator in $H_\rho$:
$$
\rho(f):=\int_G f(g)\,\rho(g)\, dg.
$$
This determines a representation of the convolution algebra $L^1(G)$ in $H_\rho$,
$$
\rho(f_1) \rho(f_2)= \rho(f_1*f_2).
$$

Consider a Borel measure $\nu$ on $\wh G$ and the direct integral of Hilbert spaces $H_\rho$
with respect to the measure $\nu$.
Consider the space $\cL(\wh G, \nu)$ of measurable functions $\Phi$ on $\wh G$
sending any  $\rho\in G$ to a Hilbert--Schmidt
operator in $H_\rho$ and satisfying the condition
$$
\int_{\wh G} \tr \bigl( \Phi(\rho)^* \Phi(\rho)\bigr)\, d\nu(\rho)<\infty.
$$

\smallskip

{\it There exists a unique measure $\mu$ on $\wh G$ (the {\it Plancherel measure}), such that for any
$f_1$, $f_2\in L^1\cap L^2(G)$ we have
$$
\la f_1,f_2\ra_{L^2(G)}=\int_{\wh G} \tr\bigl( \rho(f_2)^*\rho(f_1)\bigr)\,d\mu(\rho)
$$
and the map $f\mapsto \rho(f)$ extends to a unitary operator from $L^2(G)$ to the space
 $\cL^2(\wh G,\mu)$} (F.~I.~Mautner, I.~Segal (1950), see, e.g., \cite{Fol}).

\smallskip

{\bf 2. An example. The group $\GL(2,\R)$.} 
Let $\GL(2,\R)$ be the group of invertible real matrices 
of order $2$. 
Let $\mu\in\C$ and $\epsilon\in\Z_2$. We define the function
$x^{\mu\kos  \epsilon}$ on $\R\setminus 0$ by
$$
x^{\mu\kos  \epsilon}:=|x|^\mu \sgn(x)^\epsilon.
$$

Denote 
$\Lambda:=\C\times \Z_2\times \C\times \Z_2$.
For each element $(\mu_1,\epsilon_1;\mu_2,\epsilon_2)$ of $\Lambda$ we define a representation  
$T_{\mu,\epsilon}$  of $\GL_2(\R)$ in the space of functions 
on $\R$ by 
\begin{multline*}
T_{\mu_1,\epsilon_1;\mu_2,\epsilon_2}\begin{pmatrix}
a&b\\
c&d
\end{pmatrix} \phi(t)=\\=
\phi\Bigl(\frac{b+t d}{a+t c}\Bigr)
\cdot(a+t c)^{-1+\mu_1-\mu_2\kos  \epsilon_1-\epsilon_2} \det\begin{pmatrix}
a&b\\
c&d
\end{pmatrix}^{1/2+\mu_2\kos  \epsilon_2}
.
\end{multline*}

This formula determines the {\it principal series of representations of $\GL(2,\R)$}.
If $\mu_1-\mu_2\notin \Z$, then representations  $T_{\mu_1,\epsilon_1;\mu_2,\epsilon_2}$ and
$T_{\mu_2,\epsilon_2;\mu_1,\epsilon_1}$ are irreducible and equivalent (on representations of
$\SL(2,\R)$, see, e.g., \cite{GGV}, \cite{Vara}).

If $\mu_1=i \tau_1$, $\mu_2=i \tau_2\in i\R$, then a representation $T_{\mu_1,\epsilon_1;\mu_2,\epsilon_2}$ 
is unitary in $L^2(\R)$
(they are called representastions the {\it unitary principal series}).

Next, we define {\it representations of the discrete series}.
Let $n=1$, 2, 3, \dots. Consider the Hilbert space $H_n$ of 
holomorphic functions $\phi$ on $\C\setminus\R$ satisfying 
$$
\int_{\C\setminus\R} |\phi(z)|^2 |\Im z|^{n-1} \,d\Re z\,d\Im z<\infty.
$$
In fact, $\phi$ is a pair of holomorphic functions 
determined
  on half-planes $\Im z>0$ and $\Im z<0$.
For 
$\tau\in\R$, $\delta\in\Z_2$ we define the unitary representation $D_{n,\tau,\delta}$
of $\GL_2(\R)$ in $H_n$ by
\begin{equation*}
D_{n,\tau,\delta}\begin{pmatrix}
a&b\\
c&d
\end{pmatrix}\phi(z)= \phi\Bigl(\frac{b+z d}{a+z c}\Bigr)
(a+z c)^{-1-n}
\det \begin{pmatrix}
a&b\\
c&d
\end{pmatrix}^{1/2+n/2+i\tau\kos  \delta}
.
\end{equation*}

There exists also the complementary series of unitary representations, which does not participate in the Plancherel formula.

\smallskip

{\sc Remark.}
The expression for $D_{n,\tau,\delta}$ is contained in the family 
$T_{\mu_1,\epsilon_1;\mu_2,\epsilon_2}$, but we change the space of the representations.
\hfill $\square$

\smallskip

The Plancherel measure for $\SL(2,\R)$ was explicitly evaluated in 1952 by Harish-Chandra,
it is supported by the principal and discrete series. On the principal series the density given
by the formula (see, e.g. \cite{Vara})
\begin{align*}
d\cP=\frac 1{16\pi^3}(\tau_1-\tau_2)
\tanh \pi(\tau_1-\tau_2)/2\,d\tau_1\,d\tau_2,
\qquad \text{if $\epsilon_1-\epsilon_2=0$;}
\\
d\cP=\frac 1{16\pi^3}
(\tau_1-\tau_2)
\coth \pi(\tau_1-\tau_2)/2\,d\tau_1\,d\tau_2 \qquad \text{if $\epsilon_1-\epsilon_2=1$.}
\end{align*}

On $n$-th piece of the discrete series the measure is given by
$$
d\cP=\frac n{8\pi^3} d\tau.
$$

{\bf 3. Homogeneous spaces, etc.} The Plancherel formula for complex classical groups was obtained by
I.~M.~Gelfand and M.~A.~Naimark  \cite{GN} in 1948-50, for real semisimple groups by Harish-Chandra
in 1965 (see, e.g., \cite{Herb}, \cite{Knapp}),
there is also a formula for nilpotent groups (A.~A.~Kirillov \cite{Kir}, L.~Pukanszky \cite{Puk}).

During 1950-- early 2000s there was obtained  a big zoo of explicit spectral decompositions of $L^2$ on homogeneous spaces,
of tensor products of unitary representations,
of restrictions of unitary representations to subgroups. We present some references, which can be useful for our 
purposes \cite{Bopp}, \cite{GN}, \cite{Harinck}, \cite{Herb},  \cite{Molchanov-hyper},
\cite{Ner-beta}, \cite{Ner-Ber},  \cite{Sano}, \cite{Wallach}.
Unfortunately, texts about groups of rank $>1$ are written for experts and are heavy for exterior readers. 
See also the paper \cite{Ner-sahi}
on some spectral problems (deformations of $L^2$ on pseudo-Riemannian symmetric spaces),
which apparently are solvable but are not solved.

However, a development of the last decades seems strange. The Plancherel formula for Riemannain symmetric spaces \cite{GK}
(see, e.g., \cite{Helgason})
and Bruhat--Tits buildings \cite{Macdonald}
had a general mathematical influence (for instance to theory of special functions and to theory of integrable systems). 
Usually, Plancherel formulas are heavy results  (with impressive explicit formulas)
without further continuation even inside representation theory 
and non-commutative harmonic analysis.

\smallskip

{\bf 4. Operational calculus for $\GL(2,\R)$}, see \cite{Ner-GL2}, 2017. Denote by $\Gr_4^2$ 
the Grassmannian of all 2-dimensional linear subspaces in $\R^4$.
The natural action of the group $\GL(4,\R)$ in $\R^4$ induces the action on $\Gr_4^2$, therefore we have a unitary representation
of the group $\GL(4,\R)$ in $L^2$ on $\Gr_4^2$ (this is an irreducible representation of a degenerate principal series)
and the corresponding action of the Lie algebra $\mathfrak{gl}(4)$.

For  $g\in \GL(2,\R)$  its graph is a linear subspace in $\R^2\oplus \R^2=\R^4$.
           In this way we get an embedding 
           $$\GL(2,\R)\, \to \,\Gr_4^2.$$
           The image of the embedding is 
           an open dense subset in $\Gr_4^2$. Thus we have an identification of Hilbert spaces
           $$
           L^2\bigl(\GL(2,\R)\bigr)\simeq L^2\bigl(\Gr_4^2\bigr)
           $$
           (since natural measures on $\GL(2,\R)$ and $\Gr_4^2$ are different,
           we must multiply functions by an appropriate density  to obtain a unitary operator).
           Therefore we get a canonical action of the group $\GL(4,\R)$ 
           in $L^2\bigl(\GL(2,\R)\bigr)$. It is easily to see that the block diagonal subgroup
           $\GL(2,\R)\times \GL(2,\R)\subset \GL(4,\R)$
           acts by left and right shifts on $\GL(2,\R)$.

           We wish to evaluate the action of the Lie algebra $\mathfrak{gl}(4)$ in the Fourier-image.
           
       \smallskip    
           
Consider the space $C_0^\infty\bigl(\GL(2,\R)\bigr)$ of smooth compactly supported functions on $\GL(2,\R)$. 
For any $F\in C_0^\infty\bigl(\GL(2,\R)\bigr)$ consider the operator-valued function
$T_{\mu_1,\epsilon_1; \mu_2,\epsilon_2}(F)$ depending on 
$
(\mu_1,\epsilon_1; \mu_2,\epsilon_2)\in \Lambda
$.
We write  these operators in the form
 \begin{equation*}
 T_{\mu_1,\epsilon_1;\mu_2;\epsilon_2}(F)\phi(t)=
 \int_{-\infty}^{\infty} K(t,s|\mu_1,\epsilon_1;\mu_2,\epsilon_2)\,
 \phi(s)\,ds.
 \end{equation*}
 The kernel $K$ is smooth in $t$, $s$ and holomorphic in $\mu_1$, $\mu_2$. 

 On the other hand we have the Hilbert space $\cL^2\bigl(\wh{\GL(2,\R)}, d\cP\bigr)$.
 The norm in this Hilbert space is given by
 \begin{multline}
 \|K\|^2=\int \int_{-\infty}^\infty \int_{-\infty}^\infty \bigl|K(t,s|\mu_1,\epsilon_1;\mu_2,\epsilon_2) \bigr|^2
 dt\,ds\, d\cP(\mu) +
 \\
 + \Bigl\{ \text{summands corresponding to the discrete series}\Bigr\}
 .
 \label{eq:dcp}
 \end{multline}

 We must write 
 the action of the Lie algebra $\mathfrak{gl}(4)$.
 Denote by $e_{kl}$ the standard generators of $\mathfrak{gl}(4)$ acting in smooth compactly supported
 functions on $\GL(2,\R)$ and by $E_{kl}$ the same generators acting in the space of functions
 of variables $t$, $s$, $\mu_1$, $\epsilon_1$, $\mu_2$, $\epsilon_2$.
 The action of the subalgebra $\mathfrak{gl}(2)\oplus \mathfrak{gl}(2)$
 is clear from the definition of the Fourier transform, this Lie algebra acts by first order differential operators.
 For instance
 \begin{align*}
 e_{12}&=-b\frac\partial{\partial a}- d \frac\partial{\partial b},\qquad
 E_{12}= \frac\partial{\partial t};
 \\
 e_{43}&= b \frac\partial{\partial a}+ d \frac\partial{\partial c}, \qquad
 E_{43}=-s^2 \frac\partial{\partial s}+(-1-\mu_1+\mu_2)s.
 \end{align*}
 
 Define shift operators $V_1^+$, $V_1^-$, $V_2^+$, $V_2^-$
 by
 \begin{align}
 V_1^\pm  K(t,s|\mu_1,\epsilon_1;\mu_2,\epsilon_2)= K(t,s|\mu_1\pm 1,\epsilon_1+ 1;\mu_2,\epsilon_2);
  \label{eq:v1}
 \\
 V_2^\pm  K(t,s|\mu_1,\epsilon_1;\mu_2,\epsilon_2)= K(t,s|\mu_1,\epsilon_1;\mu_2\pm 1,\epsilon_2+ 1).
 \label{eq:v2}
 \end{align}
 To be definite, we present  formulas for two nontrivial generators $e_{kl}$ and their Fourier images $E_{kl}$:
\begin{align*}
e_{14}&=\frac \partial {\partial b}+\frac c{ad-bc},
\\
  E_{14}&=\frac{-1/2  + \mu_1}{\mu_1 - \mu_2}\, \frac\partial{\partial s}\, V_1^- 
  +\frac{-1/2  + \mu_2}{\mu_1 - \mu_2}\,\frac\partial{\partial t}\, V_2^-,
  \label{eq:1}
\\
 e_{32}&=-\Bigl(ac \frac{\partial} {\partial a}+ad\frac{\partial} {\partial b}+c^2\frac{\partial} {\partial c}+
 cd \frac{\partial} {\partial d}\Bigr)-c,
 \\
     E_{32}&=\frac{1/2 +  \mu_1 }{\mu_1 - \mu_2}\,
     \frac\partial{\partial t}\, V_1^+  +\frac{1/2 +  \mu_2 }{\mu_1 - \mu_2} \, \frac\partial{\partial s}\,
     V_2^+.
     \end{align*} 
 There is also a correspondence for operators of multiplication by functions.
   For instance,  the operator of multiplication by
 $c$ in $C_0^\infty\bigl(\GL(2,\R)\bigr)$
 corresponds to
 $$
 \frac1{\mu_1-\mu_2} \Bigl(\frac\partial {\partial t} V_1^+ + \frac\partial {\partial s} V_2^+
 \Bigr)
 $$
 in the Fourier-image.
There are similar formulas for multiplications by $a$, $b$, $d$.
 The operator of multiplication by $(ad-bc)^{-1}$ corresponds to $V_1^- V_2^-$ (the last statement is trivial).
 The operator $\frac\partial{\partial b}$
 corresponds to
   $$\frac{\mu_1-\tfrac 32}{\mu_1 - \mu_2}\, \frac\partial{\partial s}\, V_1^- 
  +\frac{\mu_2-\tfrac 32}{\mu_1 - \mu_2}\,\frac\partial{\partial t}\, V_2^-.
  $$
 There are similar formulas for other partial derivatives.

 \sm

We emphasize  that  {\bf our formulas contain shifts in imaginary directions} (the shifts in
(\ref{eq:v1})--(\ref{eq:v2}) are transversal to the contour of integration  in (\ref{eq:dcp})).

\smallskip

{\bf 5. Difference operators in imaginary direction and classical integral transforms.}
 The operators $i E_{kl}$  are symmetric in the sense of the spectral theory.
 The question about domains of self-adjointness  is open. 

There exist elements of
spectral theory of self-adjoint difference operators 
in $L^2(\R)$ of the type
\begin{equation}
L f(s)=a(s) f(s+i)+ b(s)f(s)+ c(s) f(s-i),\qquad i^2=-1,
\label{eq:i}
\end{equation}
see \cite{Ner-imaginary}, \cite{Gro}.
 Recall that several systems of classical  hypergeometric
orthogonal polynomials (Meixner-Polaszek, continuous Hahn, continuous dual Hahn, Wilson,
see, e.g. \cite{Koe})
are eigenfunctions of operators of this type. In the polynomial cases the problems are  algebraic.
The simplest nontrivial analytic  example is the operator 
$$
M f(s)=\frac 1{is} \bigl( f(s+i)-f(s-i)\bigr)
$$
in $L^2\bigl(\R_+, |\Gamma(is)|^{-2} ds\bigr)$. We define $M$ on the space of functions
$f$
holomorphic in a strip $|\Im s|<1+\delta$ and satisfying the condition
$$|f(s)|\le \exp\{-\pi |\Re s| \} |\Re s|^{-3/2-\epsilon}$$
in this strip.
The spectral decomposition of $M$
is given by the inverse Kontorovich--Lebedev integral transform. 
Recall that the direct Kontorovich--Lebedev
transform
$$
\cK f(s)=\int_0^\infty K_{is}(x) f(x) \frac {dx}{x},\qquad \text{where $K_{is}$ is  the Macdonald--Bessel function,}
$$
gives the spectral decomposition of a second order differential operator, namely 
$$D:=\bigl( x\frac d{dx}\bigr)^2-x^2, \qquad  x>0.$$
The transform $\cK$ is a unitary operator
$$L^2(\R_+,dx/x) \to L^2\bigl(\R_+, |\Gamma(is)|^{-2} ds\bigr).$$
It send
$D$ to the multiplication by $s^2$, and $\cK^{-1}$ send the difference operator $M$ to the multiplication by $2/x$.
So we get so-called {\it bispectral problem}.

Now there  is a zoo of explicit spectral decompositions
of  operators (\ref{eq:i}). The similar bispectrality appears for some other 
integral transforms:
the index hypergeometric transform (another names of this transform are: the Olevsky transform, the Jacobi transform,
the generalized
Mehler--Fock transform) \cite{Ner-index}, the  Wimp transform with Whittaker kernel \cite{Ner-imaginary},
 a continuous analog of expansion in Wilson polynomials proposed by W.~Groenevelt \cite{Gro}, etc.

This science now is a list of examples (which certainly can be extended), but there are no a priory theorems.


{\bf 6. A general problem about overalgebras.} 
 {\it  Let $G$ be a Lie group, $\mathfrak{g}$ the Lie algebra.
Let $H\subset G$ be a subgroup. Let $\sigma$ be an irreducible unitary representation
of $G$. Assume that we know an explicit spectral decomposition of restriction of $\rho$ to a subgroup
$H$. To write the action of the overalgebra $\mathfrak{g}$  in the spectral decomposition.}

\smallskip

{\sc Remarks.}
1) Above we have $G=\GL(4,\R)$, its representation $\sigma$ in $L^2$ on the Grassmannian $\Gr_4^2$, and $H=\GL(2,\R)\times \GL(2,\R)$.
The restriction problem
is equivalent to the decomposition of regular representation of $\GL(2,\R)\times \GL(2,\R)$
in $L^2\bigl( \GL(2,\R)\bigr)$. The Fourier transform is the spectral decomposition of the regular representation.

2) It is important that similar overgroups exist for all 10 series of classical real Lie groups%
\footnote{More precisely, 
an overgroup $\wt G$ exists for $G=\GL(n,\R)$, $\GL(n,\C)$, $\GL(n,\H)$, $\O(p,q)$, $\U(p,q)$,
$\Sp(p,q)$, $\Sp(2n,\R)$, $\Sp(2n,\C)$, $\O(n,\C)$, $\SO^*(2n)$  (and not for $\SL(n,\cdot)$,
$\SU(p,q)$). For instance, for $g\in \Sp(2n,\R)$ its graph is a Lagrangian subspace in $\R^{2n}\oplus \R^{2n}$,
this determines a map from $\Sp(2n,\R)$ to the Lagrangian Grassmannian with an open dense image.
We set $\wt G:=\Sp(4n,\R)$.}.
Moreover,
a decomposition of $L^2$ on any classical symmetric space%
\footnote{The  groups $G$, $M$ must be from the list of the previous footnote, $M$ must be a symmetric
subgroup in $G$.}
$G/M$ can be regarded 
as a certain restriction problem, see \cite{Ner-pseudo}. 

3) Next,  consider a tensor product $\rho_1\otimes \rho_2$ of two unitary representations of a group $G$.
Then we have the action of $G\times G$ in the tensor product, so the problem of decomposition of tensor
products can be regarded as a problem of a restriction from the group $G\times G$  to the diagonal subgroup $G$.
\hfill $\square$

\smallskip

The question under the discussion was formulated in \cite{Ner-imaginary}. Several problems of this kind
were solved \cite{Ner-imaginary}, \cite{Mol1}--\cite{Mol3}, \cite{Ner-gl}, \cite{Ner-GL2}. 
In all the cases we get differential-difference operators
including shifts in imaginary direction. Expressions also include differential operators of high order,
even for $\SL(2,\R)$-problems we usually get operators of order 2.

\smallskip

{\bf Conjecture.} {\it All problems of this kind are solvable (if we are able to write a spectral decomposition)}.

\smallskip

{\bf 6. The Gelfand-Gindikin problem,} \cite{GG}, 1977. The set $\wh H$ 
of unitary representations of a semisimple group $H$ naturally splits into 
different types (series).


{\it Let $H$ be a semisimple group, $M$ a subgroup. Consider the space $L^2(H/M)$. Usually its $H$-spectrum
contains different series. To write explicitly decomposition of $L^2$ into pieces with uniform spectrum.}

A variant of the problem: let $G$ be a  Lie group,  $H\subset G$ a semisimple subgroup, $\rho$ is a unitary representation
of $G$. Answer to the same question.

\smallskip

{\bf 7. Example: separation of series for the one-sheet hyperboloid.}
Consider the space $\R^3$ equipped with an indefinite inner product 
$$\la u,v\ra=-u_1 v_1+u_2 v_2 + u_3 v_3.$$ 
Consider the pseudo-orthogonal group
 preserving the form $\la\cdot,\cdot\ra$, denote by $\SO_0(2,1)$ its connected component.
Recall that
$\SO_0(2,1)$ is isomorphic to the quotient $\PSL(2,\R)$ of $\SL(2,\R)$ by the center $\{\pm 1 \}$.

Consider a one-sheet hyperboloid $H$ defined by $x_1^2-x_2^2-x_3^2=1$. It is 
an $\SO_0(2,1)$-homogeneous space admitting a unique (up to a scalar factor) invariant measure.
Decomposition of $L^2(H)$ into irreducible representations of $\SO_0(2,1)$ is well-known.
The spectrum is a sum of all representation of the discrete series of $\PSL(2,\R)$ and 
the integral over the whole principal series with multiplicity 2.
The separation of series was proposed by V.~F.~Molchanov \cite{Mol-im}  in 1980
(we use a modification from \cite{Ner-Mol}).

Denote by $\ov \C=\C\cup \infty$ the Riemann sphere, by $\ov \R=\R\cup \infty$ denote 
the  the real projective line, $\ov\R\subset\ov \C$.
Consider the diagonal action of $\SL(2,\R)$ on $\ov\C\times \ov\C$,
$$
(x_1,x_2)\mapsto \Bigl(\frac {b+ d x_1}{a+ c x_1}, \frac {b+ d x_2} {a+ c x_2}\Bigr).
$$
Consider the subset $H'$ in $\ov\R\times \ov\R$ consisting of points  $x_1$, $x_2$
such that $x_1\ne x_2$. 
 It is easy to verify that $H'$ is an orbit of $\SL(2,\R)$,
it is  equivalent to the hyperboloid $H$ as a homogeneous space%
\footnote{Two families of lines on the hyperboloid correspond to two families of
 lines $x_1=\mathrm{const}$ and $x_2=\mathrm{const}$ on $\ov \R\times \ov \R$.}. It is easy to verify that
the invariant measure on $H'$ is given by the formula
$$
d\nu(x_1,x_2)=
|x_1-x_2|^{-2} \,dx_1 \,dx_2.
$$
We identify the space $L^2(H',d\nu)$ with the standard $L^2(\R\times \R)$  by the {\it unitary} operator
$$
J f(x_1,x_2)= f(x_1,x_2) (x_1-x_2)^{-1}.
$$
Now our representation in $L^2(H)$ transforms to the following 
unitary representation in 
the standard $L^2(\R^2)$: 
\begin{equation}
Q\begin{pmatrix}a&b\\c & d \end{pmatrix}
f(x_1,x_2)= f\Bigl(\frac {b+ d x_1}{a+ c x_1}, \frac {b+ d x_2} {a+ c x_2}\Bigr)
(a+ c x_1)^{-1} (a+c x_2)^{-1}.
\label{eq:}
\end{equation}

Next, consider a unitary reprepresentation of $\SL(2,\R)$ in $L^2(\R)$ given by
$$
T\begin{pmatrix}a&b\\c & d \end{pmatrix}
f(x)= f\Bigl(\frac{b+xd}{a+xc}\Bigr) (a+xc)^{-1}.
$$
Obviously, we have
$$
Q=T\otimes T.
$$
The representation $T$ is contained in the unitary principal series and it is a unique reducible element
of this series (see, e.g., \cite{GGV}). 

Denote by $\Pi_\pm$ the upper and lower half-planes in $\ov\C$.
The Hardy space $H^2(\Pi_+)$ consists of functions $F_+$ holomorphic in $\Pi_+$
that can be represented in the form
$$
F_+(x)=\int_0^\infty \phi(t) e^{itx}\,dt, \qquad\text{where $\phi(t)\in L^2(\R_+)$}.
$$
Obviously, $F$ is well-defined also on $\R$ and is contained in $L^2$.
The space $H^2(\Pi_-)$ consists of functions $F_-$ holomorphic in $\Pi_-$
of the form
$$
F_-(x)=\int_{-\infty}^0 \phi(t) e^{itx}\,dt, \qquad\text{where $\phi(-t)\in L^2(\R_+)$}.
$$
Evidently,
$$
L^2(\R)=H^2(\Pi_+)\oplus H^2(\Pi_+).
$$
It can be shown that the subspaces $H^2(\Pi_\pm)\subset L^2(\R)$ are invariant with respect to operators $T(\cdot)$,
and therefore $T$ splits into two summands $T_+\oplus T_-$ (one of them has a highest weight, another a lowest weight).
Hence 
$$Q=(T_+\oplus T_-)\otimes (T_+\oplus T_-)$$
splits into 4 summands.
It can be shown that this is the desired decomposition:


--- the space $H^2(\Pi_+)\otimes H^2(\Pi_+)$ consists of functions in $L^2(\R^2)$ 
continued holomorphically to the domain $\Pi_+\times \Pi_+$;
the representation $T_+\otimes T_+$ in $H^2(\Pi_\pm)\subset L^2(\R)$
 is a direct sum of all highest weight representations
of representation of $\PSL(2,\R)$;

--- $T_-\otimes T_-$  is a direct sum of all lowest weight representations; 

--- in $T_+\oplus T_-$ we have the direct integral of all representations of the principal series
(and the same integral in $T_-\otimes T_+$).


{\sc Remark.}  S.~G.~Gindikin \cite{Gindikin} used a similar argument
(restriction from a reducible representation of an overgroup) for
multi-dimensional hyperboloids. \hfill $\square$.

\sm

{\bf 8. Splitting off the complementary series,} see \cite{NO}.
Consider the pseudo-orthogonal group $\O(1,q)$ consisting of operators preserving
the following indefinite inner product
in $\R^{1+q}$,
$$
\la x, y\ra= -x_0y_0 + x_1y_1+\dots+ x_q y_q. 
$$
We write  elements of this group as block $(1+q)\times (1+q)$ matrices $g=\begin{pmatrix}a&b\\c&d \end{pmatrix}$.
Denote by $\SO_0(1,q)$ its connected component, it consists of matrices satisfying 
two additional conditions $\det g=+1$, $a>0$.
Denote by $S^{q-1}$ the unit sphere in $\R^n$. The group $\O(1,q)$ acts
on $S^{q-1}$ by conformal transformations $x\mapsto (a+xc)^{-1} (b+xd)$
(they preserve the sphere), the coefficient of a dilatation equals to
$(a+xc)^{-1}$.

For  $\lambda\in \C$ we define a representation $T_\lambda= T_\lambda^q$ of $\SO_0(1,q)$ 
in a space of functions on $S^{q-1}$ by
$$
T_\lambda \begin{pmatrix}a&b\\c&d \end{pmatrix} f(x)= (a+xc)^{-(q-1)/2+\lambda} f\bigl((a+xc)^{-1} (b+xd)\bigr).
$$
If $\lambda=i\sigma\in i \R$, then our representation is unitary in $L^2(S^{q-1})$, 
in this case $T_{i\sigma}$ is called a representation of the {\it unitary spherical principal series},
representations
$T_{i\sigma}$ and $T_{-i\sigma}$ are equivalent
(on these representations see e.g. \cite{Vilenkin}). If $0<s<(q-1)/2$, then $T_s$ is unitary
in the Hilbert space $H_s$ with the 
the inner product
$$
\la f_1,f_2\ra_s=\int_{S^{q-1}} \int_{S^{q-1}} \frac{f_1(x_1)\,\ov{f_2(x_2)}\,dx_1\,dx_2 }
{\|x_1-x_2|^{(q-1)/2-s}}.
$$
More precisely, $\la,\cdot,\cdot \ra$ determines a positive definite Hermitian 
form on the space $C^\infty(S^{q-1})$  (this is not obvious), we get
a pre-Hilbert space and consider its 
completion $H_s$. Such representations form the {\it spherical complementary series}.
The spaces $H_s$ are Sobolev spaces%
\footnote{In the standard notation, $H_s$ is the Sobolev space $H^{-s,2}(S^{q-1})$. Notice that 
Sobolev spaces $H^{\sigma,2}(\cdot)$ are Hilbert spaces but inner product are  defined not canonically.
In our case the inner products are uniquely determined from the $\SO_0(1,q)$-invariance.
For semisimple groups of rank $>1$ complementary series are realized in functional Hilbert spaces that are not
Sobolev spaces.}.

Consider a restrictions of $T_{i\sigma}$ to the subgroup $\SO_0(1,q-1)$.
The group $\SO_0(1,q-1)$ has the following orbits on $S^{q-1}$: the equator $Eq=S^{q-2}$
defined by the equation $x_q=0$, the upper hemisphere $H_+$
and the lower hemisphere $H_-$. The equator has zero measure and can be forgotten.
Therefore 
$$L^2(S^{q-1})=L^2(H_+)\oplus L^2(H_-).$$
On the other hand,
hemispheres as homogeneous spaces are equivalent to $\SO_0(1,q-1)/\SO(q-1)$,
i.e. to the $(q-1)$-dimensional  Lobachevsky space. The decomposition of $L^2$ is a classical problem,
in each summand $L^2(H_\pm)$ we get a multiplicity-free direct integral over the whole
spherical principal series.

\sm

The restriction of a representation $T_s$ of the complementary series
is more interesting, it contains several summands of the complementary series
and is equivalent to
\begin{equation}
\bigoplus\nolimits_{k: \, s-k>1/2} T^{q-1}_{s-k}\,\, \bigoplus\,\, L^2(H_+)\,\,\bigoplus\,\, L^2(H_-).
\label{eq:decomposition}
\end{equation}
This spectrum was obtained by Ch.~Boyer (1973), our purpose is to visualize summands of the complementary series.

According  {\it trace theorems} Sobolev spaces of negative order can contain
distributions supported by submanifolds.
Denote by $\delta_{Eq}$ the delta-function of the equator,
$\delta_{Eq}:=\delta(x_q)$.
Let $\phi$ be a smooth function $\phi$ on $Eq$.
$$
\|\phi \delta_{Eq}\|^2_{s}=
\la \phi \delta_{Eq}, \phi \delta_{Eq}\ra_s=
\int_{S^{q-2}} \int_{S^{q-2}} \frac{\phi(y_1)\,\ov{\phi(y_2)}\,dy_1\,dy_2 }
{\|y_1-y_2|^{-(q-1)/2+s}}.
$$
If $s>1/2$ the integral converges and
 $\phi \delta_{Eq}\in H_s$. The representation of $\SO_0(1,q)$
 in the space of such functions is $T_s^{q-1}$.
 
 Denote by $\frac\partial {\partial n} \delta_{Eq}:=\delta'(x_q)$ the  derivative 
 of $\delta_{Eq}$ in the normal direction. Similar arguments show that for $s>3/2$ and smooth $\psi$ we have
 $\psi \frac\partial {\partial n} \delta_{Eq}\in H_s$. The space  of functions of the form
 $$\phi \delta_{Eq}+\psi \frac\partial {\partial n} \delta_{Eq}$$
 again is invariant. It contains the subspace $T_s^{q-1}$ and we get
 the representation $T_{s+1}^{q-1}$ in the quotient. Since our representation is unitary, 
 $T_{s+1}^{q-1}$ must be direct summand. Etc.

Next, we consider the operator $J:H_s\mapsto L^2(S^{q-1})$
given by
$$Jf(x)=|x_q|^{(q-1)/2-s}f(x).$$
It  intertwines restrictions of $T_s$ and $T_0$, the kernel of $J$ consists of 
distributions supported by $Eq$ and the image is dense%
\footnote{More precisely, we consider this operator as an operator on smooth functions 
compactly supported outside $Eq$, take the closure $\Gamma$ of its graph in $H_s\oplus L^2$,
and examine projection operators $\Gamma\to H_s$, $\Gamma\to L^2$.}. This gives us 
(\ref{eq:decomposition}).

\sm

{\bf 9 The modern status of the problem.}
We mention the following works:

a) G.~I.~Olshanski \cite{Olsh} (1990) proposed a way to split off highest weight and lowest weight representations.

b) The author in \cite{Ner-old} (1986)  proposed a way to split off complementary series 
(see proofs and further examples  in \cite{NO},
the paper \cite{Ner-separation} contains an example with separation
 of direct integrals of different complementary series).

c) S.~G.~Gindikin \cite{Gindikin} (1993) and  V.~F.~Molchanov \cite{Mol-separ} (1998)  obtained
 a separation of spectra for multi-dimensional hyperboloids.

\smallskip

These old works had continuations, in particular were many further works 
with splitting off highest weight representations
(for more references, see \cite{Ner-Upq}).

The recent paper  \cite{Ner-Upq}  (2017) contains formulas for
projection operators separating spectrum  for $L^2$ on pseudo-unitary groups $\U(p,q)$.
In this case we can consider separation into series (if we fix the number $r$ of continuous parameters
of a representation, $r\le \min(p,q)$),  subsubseries (if we fix all discrete parameters of 
a representation) and intermediate subseries. All these  question are solvable.
The solution was obtained by a summation of all characters corresponding to a given type of spectrum, certainly
this way must be available for all semisimple
Lie groups. 

In \cite{Ner-separationGL} the problem was solved for $L^2$ on pseudo-Riemannian symmetric spaces
$\GL(n,\C)/\GL(n,\R)$. The calculation is based on an explicit summation of spherical distributions.
Apparently, this can be extended to all symmetric spaces of the form $G_\C/G_\R$, where
$G_\C$ is a complex semisimple Lie group and $G_\R$ is a real form of $G_\C$
(on Plancherel formulas for such spaces, see \cite{Bopp}, \cite{Harinck}, \cite{Sano}).

For  arbitrary
semisimple symmetric spaces the problem does not seem well-formulated,
see a discussion of multidimensional 
 hyperboloids  in \cite{Mol-separ}.

 \tt
 
 \noindent
 Math. Dept., University of Vienna; \\
 Institute for Theoretical and Experimental Physics (Moscow); \\
 MechMath Dept., Moscow State University;\\
 Institute for Information Transmission Problems.\\
 URL: http://mat.univie.ac.at/$\sim$neretin/

\end{document}